\documentclass[11pt]{amsart}

\usepackage{amssymb,amsbsy,amsmath,amsfonts,amssymb,amscd,amsthm}

\newcommand{\Cc}{\mathbb{C}}  % field of complex numbers
\newcommand{\Pp}{\mathbb{P}}
\newcommand{\Rr}{\mathbb{R}}
\newcommand{\Qq}{\mathbb{Q}}
\newcommand{\Zz}{\mathbb{Z}}
\newcommand{\defi}[1]{\emph{#1}}

\renewcommand{\epsilon}{\varepsilon}

\newcommand{\Baff}{{\mathcal{B}_\mathit{\!aff}}}
\newcommand{\Binf}{{\mathcal{B}_\infty}}
\newcommand{\B}{{\mathcal{B}}}

\newcommand{\Aut}{\mathop{\mathrm{Aut}}\nolimits}

\newcommand{\id}{\mathop{\mathrm{id}}\nolimits}

\newcommand{\mm}{{\mathfrak{m}}}

\newcommand{\CC}{{\mathcal{C}}}

{
\newtheorem*{theorem*}{Theorem}  % theorem without number
\newtheorem*{theoremA}{Theorem A}  % theorem without number
\newtheorem*{theoremB}{Theorem B}  % theorem without number
\newtheorem{lemma}{Lemma}        % lemma with number
}

\begin{document}

\title{%
Non reality and non connectivity of complex polynomials
}
\author{%
Arnaud Bodin
}
\address{%
Laboratoire Agat,
UFR de Math\'ematiques,
Universit\'e Lille 1, 59655 Villeneuve d'Ascq Cedex, France.
E-mail: Arnaud.Bodin@agat.univ-lille1.fr
}
\begin{abstract}{%
Using the same method we provide negative answers to the following
questions: Is it possible to find real equations for complex
polynomials in two variables up to topological equivalence (Lee 
Rudolph)? Can two topologically equivalent polynomials be 
connected by a continuous family of topologically equivalent 
polynomials? }
\end{abstract}

\maketitle

\section{Introduction}
 Two polynomials $f,g \in \Cc[x,y]$ are \defi{topologically
equivalent}, and we will denote $f\approx g$, if there exist
homeomorphisms  $\Phi : \Cc^2 \longrightarrow \Cc^2$ and $\Psi 
:\Cc \longrightarrow \Cc$
 such that $g\circ \Phi = \Psi \circ f$.
They are \defi{algebraically equivalent}, and we will denote
$f\sim g$, if we have  $\Phi \in \Aut \Cc^2$ and 
$\Psi=\id$.

It is always possible to find real equations for germs of plane
curves up to topological equivalence. In fact the proof is as
follows: the topological type of a germ of plane curve $(C,0)$ is
determined by the characteristic pairs of the Puiseux expansions
of the irreducible branches and by the intersection multiplicities
between these branches. Then we can choose the coefficients of the
Puiseux expansions in $\Rr$ (even in $\Zz$). Now it is possible
(see \cite{EN}, appendix to chapter 1) to find a polynomial $f \in
\Rr[x,y]$ (even in $\Zz[x,y]$) such that the germ $(f=0,0)$ is
equivalent to the germ $(C,0)$.

This property has been widely used by N. A'Campo and others (see
\cite{AC} for example) in the theory of divides. Lee Rudolph asked
the question whether it is true for polynomials \cite{Ru}. We 
give a negative answer:
\begin{theoremA}
Up to topological equivalence it is not always possible to find
real equations for complex polynomials.
\end{theoremA}

\section{}
 We now deal with another problem. In \cite{Bo3} we proved 
that a family of polynomials with some constant numerical data 
are all topologically equivalent. More precisely for a polynomial 
let $\mm = (\mu, \#\Baff, \lambda, \#\Binf, \#\B)$ be the 
multi-integer respectively composed of the affine Milnor number, 
the number of affine critical values, the Milnor number at 
infinity, the number of critical values at infinity, the number 
of critical values (with $\B = \Baff \cup \Binf$). Then we have a 
global version of the L\^{e}-Ramanujam $\mu$-constant theorem:

\begin{theorem*}[\cite{Bo3}]
 Let $(f_t)_{t\in[0,1]}$ be a family of complex
polynomials in two variables  whose coefficients are polynomials
in $t$. Suppose that the multi-integer $\mm(t)$ and the degree
$\deg f_t$ do not depend on $t \in[0,1]$.  Then the polynomials
$f_0$ and $f_1$ are topologically equivalent.
\end{theorem*}

It is true that two topologically equivalent polynomials have the
same multi-integers $\mm$. A natural question is: can two
topologically equivalent polynomials  be connected by a continuous
family of topologically equivalent polynomials ?

\begin{theoremB}
There exist two topologically equivalent polynomials $f_0, f_1$
that cannot be connected by a family of equivalent polynomials.
That means that for each continuous family $(f_t)_{t\in[0,1]}$ 
there exists a $\tau \in]0,1[$ such that $f_\tau$ is not 
topologically equivalent to $f_0$.
\end{theoremB}

It can be noticed that the answer is positive for algebraic
equivalence. Two algebraically equivalent polynomials can be
connected by algebraically equivalent polynomials since $\Aut \Cc^2$
is connected by Jung's theorem.

Such kinds of problems have been studied by V. Kharlamov and V.
Kulikov in \cite{KK} for cuspidal projective curves. They give two
complex conjugate projective curves that are not isotopic. The
example with lowest degree has degree $825$. In \cite{ACA}, E.
Artal, J. Carmona and L. Cogolludo give examples of projective
curves $C,C'$ of degree $6$ that have conjugate equations in
$\Qq(\sqrt 2)$ but the pairs $(\Pp^2,C)$ and $(\Pp^2,C')$ are not
homeomorphic by an orientation-preserving homeomorphism.

\section{} 
The method used in this note is based on the relationship 
between topological and algebraic equivalence: we set a family 
$(f_s)_{s\in \Cc}$ of polynomials such that $(f_s=0)$ is a line 
arrangement in $\Cc^2$. One of the line depends on a parameter 
$s\in \Cc$. There are enough lines in order that each polynomial 
is algebraically essentially unique. Moreover every polynomial 
topologically equivalent to $f_s$ is algebraic equivalent to a 
$f_{s'}$, where $s'$ may be different from $s$.

For generic parameters the polynomials are topologically
equivalent all together and the function $f_s$ is a Morse 
function on $\Cc^2\setminus f_s^{-1}(0)$. We choose our 
counter-examples with non-generic parameters, for such an example 
$f_k$ is not a Morse function on $\Cc^2\setminus f_k^{-1}(0)$. 
The fact that non-generic parameters are finite enables us to 
prove the requested properties.

\section{Non reality}

Let
$$f_s(x,y) = xy(x-y)(y-1)(x-sy).$$
Let $k, \bar k$ be the roots of $s^2-s+1$.
\begin{theoremA}
There does not exist a polynomial $g$ with real coefficients such
that $g \approx f_k$.
\end{theoremA}
Let $\CC = \{0,1,k,\bar k\}$. Then for $s \in \Cc \setminus \CC$, 
$f_s$ verifies $\mu = 14$, $\# \Baff = 3$ and $\Binf = 
\varnothing$. By the connectivity of $\Cc \setminus \CC$ and the 
global version of the $\mu$-constant theorem, two polynomials 
$f_s$ and $f_{s'}$, with $s,s'\notin \CC$, are topologically 
equivalent.

The polynomials $f_k$ and $f_{\bar k}$ verify $\mu = 14$, but
$\#\Baff = 2$.  Then such a polynomial is not topologically
equivalent to a generic one $f_s$, $s\notin \CC$. In fact for $s
\notin \CC$ there are two non-zero critical fibers with one double
point for each one. For $s=k$ or $s= \bar k$, there is only one
non-zero critical fiber with an ordinary cusp.

\begin{lemma}
\label{lem:alg}
 Let $s,s' \in \Cc$. The polynomials $f_s$ and $f_{s'}$
are algebraically equivalent if and only if $s = s'$ or $s =
1-s'$.
\end{lemma}
In particular the polynomials $f_k$ and $f_{\bar k}$ are
algebraically equivalent.

\begin{proof}
Let us suppose that $f_s$ and $f_{s'}$ are algebraically
equivalent. Then we can suppose that there exists $\Phi \in \Aut
\Cc^2$ such that $f_{s'} = f_s \circ \Phi$. Such a $\Phi$ must
send the lines $(x=0), (y=0)$ to two lines, then $\Phi$ is linear:
$\Phi(x,y) = (ax+by,cx+dy)$. A calculus proves that
$\Phi(x,y)=(x,y)$ or $\Phi(x,y)=(y-x,y)$ that is to say $s=s'$ or
$s=1-s'$.
\end{proof}

\begin{lemma}
\label{lem:top}
 Fix $s\in\Cc$ and let $f$ be a polynomial such that $f
\approx f_s$. There exists $s'$ such that $f \sim f_{s'}$.
\end{lemma}
Then lemma \ref{lem:alg} implies that there are only two choices
for $s'$, but $s'$ can be different from $s$.

\begin{proof}
The curve $f_s^{-1}(0)$ contains the simply connected curve
$xy(x-y)(x-sy)$, then the curve $f^{-1}(0)$ contains also a simply
connected curve (with $4$ components), by the generalization of
Za\u \i denberg-Lin theorem (see \cite{Bo2}) this simply connected
curve is algebraically equivalent to $xy(x-y)(x-s'y)$. Then the
polynomial $f$ is algebraically equivalent to
$xy(x-y)(x-s'y)P(x,y)$.
The curve $C$ defined by $(P=0)$ is homeomorphic to $\Cc$ and admits
a polynomial parameterization $(\alpha(t),\beta(t))$ with $\alpha,\beta \in \Cc[t]$.
Since $C$ does not intersect the axe $(y=0)$, $\beta$ is a constant polynomial;
 since  $C$ intersects the axe $(x=0)$ at one point  $\alpha$ is monomial.
An equation of $P$ is now $P(x,y) = y^n-\lambda$. By the irreducibility of $C$ and
up to an homothety we get $P(x,y) = y - 1$. That
is to say $f$ is algebraically equivalent to $f_{s'}$.
\end{proof}

\section{}
Let $g \in \Cc[x,y]$, if $g(x,y) = \sum{a_{i,j}x^iy^j}$ then we
denote by $\bar g$ the polynomial defined by $\bar g(x,y)=
\sum{\bar a_{i,j}x^iy^j}$. Then  $g=\bar g$ if and only if all 
the coefficients of $g$ are real.

We prove theorem A. Let suppose that there exists a polynomial 
$g$ such that $g=\bar g$  and  $g \approx f_k$. There exists 
$s \in \Cc$ such that $g \sim f_s$. Since $f_k$  has only two 
critical values,  $g$ and $f_s$ have two critical values. Then 
$s=k$ or $s=\bar k$ ($s=0$ or $s=1$ gives a polynomial with 
non-isolated singularities). As $f_k \sim f_{\bar k}$ we can 
choose $s=k$. As a consequence we have $\Phi \in \Aut \Cc^2$ such 
that $ g = f_k \circ \Phi.$

 Let $\Phi$ be $\Phi = (p,q).$ Then $g 
= pq(p-q)(q-1)(p-kq)$.  As $g = \bar g$ we have :
$$ \big\{p,q,p-q,q-1,p-kq \big\} =  \big\{\bar p,\bar q,
\bar p-\bar q,\bar q-1,\bar p-\bar k \bar q \big\}.$$ 
Moreover by 
the configuration of the lines we have that $q-1 = \bar q -1$. So 
$q = \bar q$. Hence $q \in \Rr[x,y]$. So
$$ \big\{p,p-q,p-kq \big\} =  \big\{\bar p,
\bar p-\bar q,\bar p-\bar k \bar q \big\}.$$ Let suppose that
$p\not= \bar p$. Then $p = \bar p - q$ or $p = \bar p - \bar k q$.
So $p -\bar p$ equals $-q$ or $-\bar kq$. But $p-\bar p$ has
coefficients in $i\Rr$, which is not the case of $q \in \Rr[x,y]$
nor of $\bar k q$. Then $p = \bar p$. We have proved that $\Phi =
(p,q)$ has real coefficients. From $g = f_k \circ \Phi$ we get
$\bar g = \bar f_k \circ \bar \Phi$. So $g = f_{\bar k} \circ
\Phi$. On the one hand $f_k = g \circ \Phi^{-1}$ and on the
other hand $f_{\bar k} = g \circ \Phi^{-1}$. So $f_k=f_{\bar k}$,
then $k = \bar k$ which is false. It ends the proof.

We could have end in the following way: $\Phi = (p,q)$ is in $\Aut
\Cc^2$ with real coefficients, then $\Phi$, considered as a real
map, is in $\Aut \Rr^2$ (see \cite[Theorem 2.1]{Ba} for example).
Then $f_k = g \circ \Phi^{-1}$ with $g, \Phi^{-1}$ with real
coefficients, then $f_k$ has real coefficients which provides the
contradiction.

\section{Non connectivity }

Let
$$f_s(x,y) = xy(y-1)(x+y-1)(x-sy).$$
Let $\CC$ be the roots of
$$s(s-1)(s+1)(256s^4+736s^3+825s^2+736s+256)(256s^4+448s^3+789s^2+448s+256).$$
Then for $s \in \Cc \setminus \CC$, $f_s$ verifies $\mu = 14$, 
$\# \Baff = 4$ and $\Binf = \varnothing$. For $s,s'\notin \CC$, 
$f_s$ and $f_{s'}$ are topologically equivalent. The roots of 
$256s^4+448s^3+789s^2+448s+256$ are of the form $\left\lbrace k, 
\bar k,  1/k , 1/{\bar k}\right\rbrace.$ The polynomials $f_k$ 
and $f_{\bar k}$ verify $\mu = 14$, but $\#\Baff = 3$. Then such 
a polynomial is not topologically equivalent to a generic one 
$f_s$, $s\notin \CC$

\begin{theoremB}
The polynomials  $f_k$ and $f_{\bar k}$ are topologically
equivalent  and it is not possible to find a continuous family
$(g_t)_{t\in[0,1]}$ such that $g_0 = f_k$, $g_1=f_{\bar k}$ and
$g_t \approx f_k$ for all $t\in [0,1]$.
\end{theoremB}

The polynomials $f_k$ and $f_{\bar k}$ are topologically
equivalent since we have the formula $ f_{\bar k}(\bar x, 
\bar y) = \overline{f_k(x,y)}$.

The two following lemmas are similar to lemmas \ref{lem:alg} and
\ref{lem:top}.
\begin{lemma}
The polynomials $f_s$ and $f_{s'}$ are  algebraically equivalent
if and only if $s =s'$ or $s = 1/s'$.
\end{lemma}

\begin{lemma}
\label{lem:topbis}
 Fix $s$ and let $f$ be a polynomial such that $f
\approx f_s$. Then there exists $s'$ such that $f \sim f_{s'}$.
\end{lemma}

\section{}
We now prove theorem B. Let us suppose that such a family $(g_t)$ 
does exist. Then by lemma \ref{lem:topbis} for each $t\in[0,1]$ 
there exists  $s(t)\in \Cc$ such that $g_t$ is algebraically 
equivalent to $f_{s(t)}$ (in fact there are two choices for 
$s(t)$). We can suppose that there exists $\Phi_t \in \Aut \Cc^2$ 
such that $f_{s(t)} = g_t \circ \Phi_t$.

We now prove that the map $t \mapsto \Phi_t$ can be chosen
continuous, that is to say the coefficients of the defining
polynomials are continuous functions of $t$. We write $g_t =
A_tB_tG_t$ such  that $A_0(x,y) = x$, $B_0(x,y)=y$ and the maps $t
\mapsto A_t$, $t \mapsto B_t$ are continuous. So the automorphism
$\Phi_t^{-1}$ is defined by
$$\Phi_t^{-1}(x,y) = \big( A_t(x,y), B_t(x,y) \big).$$
By the inverse local theorem with parameter $t$, we have that
$t\mapsto \Phi_t$ is a continuous function. Then the map $t
\mapsto f_{s(t)}$ is a continuous function, as the composition of
two continuous functions. As $s(t)$ is a coefficient of the
polynomial $f_{s(t)}$, the map $t \mapsto s(t)$ is a continuous
function.

As a conclusion we have a map $t \mapsto s(t)$ which is continuous
and such that $s(0) = k$ and $s(1) = \bar k$. It implies that
there exists $\tau \in ]0,1[$ such that $s(\tau) \notin\CC$. On the one hand
$g_\tau$ is algebraically, hence topologically, equivalent to
$f_{s(\tau)}$; on the other hand $g_\tau$ is topologically
equivalent to $f_k$ (by hypothesis). As $s(\tau) \notin \CC$,
$f_{s(\tau)}$ and $f_k$ are not topologically equivalent (because
$\# \Baff$ are different), it provides a contradiction.

\section{}
The calculus have been done with the help of \textsc{Singular},
\cite{sing}, and especially with author's library \texttt{critic}
described in \cite{Bo4}.

This research has been done at the \emph{Centre de Recerca Matem\`atica} of Barcelona and was
supported by a Marie Curie Individual
Fellowship of the European Community (HPMF-CT-2001-01246).

%%%  Bibliography  %%%
%%%%%%%%%%%%%%%%%%%%%%%

%
\end{document}